\documentstyle{article}

\begin{document}

\newtheorem{th-ud}{Theorem}
\newtheorem{lem-ud}{Lemma}
\newcommand{\qed}{\par\hspace{11cm}\rule{5pt}{5pt}\vspace{1em}}

\title{
\hfill\raisebox{0.0ex}[1.ex][3ex]{}
              \newline
      {\bf
Almost periodic functions in finite-dimensional
space with the spectrum in a cone}
\author {Favorov S., Udodova O.${}^1$\footnotemark[0]}}
\date{}
\maketitle
\footnotetext[1]{The work was supported by the INTAS-grant ü99-00089}
\centerline
{\begin{tabular}{c}
Kharkiv National University\\
Department of Mechanics and Mathematics\\
md. Svobody 4, Kharkiv 61077, Ukraine\\
e-mail: favorov@ilt.kharkov.ua, \\
Olga.I.Udodova@univer.kharkov.ua\\
\end{tabular}}
\begin{abstract}
{ We  prove that an almost periodic function in finite-dimensional
space extends to a holomorphic bounded function in a tube domain
with a cone in  the
base if and only if the spectrum belongs to
the conjugate cone. We also prove that an almost periodic function
in finite-dimensional space has the bounded spectrum if and only if
it extends to an entire function of exponential type. \\
{ \it 2000 Mathematics Subject Classification 42A75, 30B50.}}
\end{abstract}

\bigskip

A continuous function $f(z)$ on a strip
$$
S_{a,b}=\{z=x+iy: x \in{\bf R}, \ a\le y\le b \},
\ -\infty\le a\le b\le +\infty.$$
is called {\it almost periodic by Bohr} on this strip, if
for any $\varepsilon>0$ there exists $l=l(\varepsilon)$ such that
every interval of the real axis of length $l$ contains a number
$ \tau $ ($\varepsilon$-almost period for $f(z)$)
with the property
\begin{equation}
\label{p}
\sup_{z\in S_{a, b}}|f(z+\tau)-f(z)|<\varepsilon.
\end{equation}

In particular, when $a=b=0$ we obtain the class of almost periodic functions
on the real axis.

To each almost periodic function $f(z)$ assign the Fourier series
$$
\sum\limits_{n=0}^{\infty}a_n(y)
e^{i\lambda_n x}, \ \lambda_n\in{\bf R},$$ where
$a_n(y)$ are continuous functions of the variable $y\in [a,b].$

In the case $a=b=0$ all exponents $\lambda_n$ are nonnegative
if and only if the function $f(x)$ extends to the upper half-plane
as a holomorphic bounded almost periodic function;
the set of all exponents $\lambda_n$ is
bounded if and only if $f(x)$
extends to the plane ${\bf C}$ as an entire function of the
exponential type
$\sigma=\sup\limits_n|\lambda_n|$, which is
almost periodic in every
horizontal strip of finite width (see \cite {Borr1}, \cite {Levitan}).

A number of works connected with almost
periodic functions of many variables on a {\it tube set}
appeared  recently (see \cite{Fav},
\cite{Ron0}-\cite {Ron2}).
Recall that the set $T_K\subset {{\bf C} ^m} $ is a {\it tube set}
if
$$T_{K}=\{z=x+iy: x\in {{\bf R} ^m}, y\in K \}, $$
where $K\subset {\bf R}^m$ is {\it the base} of the
tube set. \\
{\bf Definition.}
(See \cite{Ron0}, \cite{Ron2}).
{\it
\label {®¯à3}
A continuous function $f(z) $, $z\in T_K $ is called
almost periodic by Bohr on $T_K $,
if for any $\varepsilon>0$ there exists $l=l (\varepsilon) $, such that
every $m$-dimensional cube on ${\bf R}^m$ with the side $l$ contains
at least one point
$\tau$ ($(T_K,\varepsilon)$-almost period from $f (z) $)
with the property
\begin {equation}
\label {p1}
\sup\limits_{z\in T_K}|f(z+\tau)-f(z)|<\varepsilon.
\end {equation}}

Let $f,\ g$ be locally integrable functions on every
real plane
$$\{z=x+iy_0: x\in {{\bf R}^m} \}, \ y_0\in K.$$
{\bf Definition.} {\it
Stepanoff distance of the order $p\ge 1$ between functions $f$
and $g $ is the value
$$
S_{p, T_{K}}(f,g)=\sup\limits_{z\in T_{K}}\left(
\int\limits_{[0,1]^m}^{}|f(z+u)-g(z+u)|^p du\right)^
{ \frac 1p}.
$$}

Using this definition,
we can extend the concept of almost periodic functions by Stepanoff
on a strip (see~\cite {Levitan} p.197) to almost periodic functions on a
tube set: \\
{\bf Definition.} {\it
A function $f(z)$, $z\in T_K $, is called almost
periodic by Stepanoff on $T_K $,
if for any $ \varepsilon>0$ there exists
 $l=l (\varepsilon) $ such that
every $m$-dimensional cube with the side $l $ contains at least one
$\tau$ ($(T_K,\varepsilon,p)$-almost period by Stepanoff of the
function $f(z)$) with the property
\begin{equation}
\label{st}
S_{p, T_{K}}(f(z),f(z+\tau))<\varepsilon.
\end{equation}}

The Fourier series  for an
almost periodic (by Bohr or by Stepanoff) function $f(z) $ on a set
$T_K$ is the series
\begin{equation}
\label{ff}
\sum_{\lambda \in{\bf R}^m}^{}
a(\lambda,y)e^{i\langle x,\lambda\rangle},
\end {equation}
where $ \langle x, \lambda\rangle $ is the scalar product on $ {\bf R} ^m $,
and
\begin{equation}
\label{K}
a(\lambda,y)=\lim\limits_{N \to \infty}
\left(
\frac {1}{2N}
\right)^m
\int\limits_{[-N,N]^m}^{}f(x+x'+iy)e^{-i\langle x+x',\lambda\rangle} dx;
\end {equation}
this limit exists uniformly in the parameter $x'\in{\bf R}^m$ and
 does not depend on this parameter (see~\cite{Ron0}, \cite{Ron1}).

A set of all vectors $\lambda\in{{\bf R}^m}$ such that $a(\lambda,y)\not\equiv 0$
is called the spectrum of $f(z)$ and is denoted by $sp\,f $;
 this set is at most countable, therefore
the series~(\ref {ff}) can be written in the form
$$
\sum_{n=0}^{\infty}
a_n(y)e^{i\langle x,\lambda_n\rangle}.$$

Note that partial sums of the series~(\ref {ff}),
generally speaking, do not converge to the
 function $f(z)$. However the Bochner-Feyer sums
\footnotemark {}
\footnotetext {For $n=1$ see \cite{Levitan}, for $n>1$ see \cite{Ron1}.}
$$
\sigma_q (z) = \sum\limits _ {n=0} ^ {q-1}
k^q_n a_n (y) e ^ {i\langle x, \lambda_n\rangle},\ 0\le k^q_n<1, \ k^q_n
\rightarrow 1 \ as\ q\rightarrow\infty$$
converge
to the function $f(z)$
uniformly for almost periodic functions by Bohr
and with respect to the metric $S_{p,T_K}$
for almost periodic functions by Stepanoff;
in particular, if two functions have the same Fourier series,
then these functions
coincide identically.
For holomorphic almost
periodic functions the series~(\ref {ff}) can be written in the form
\begin{equation}
\label{f+}
\sum_{n=0}^{\infty}
a_n e^{-\langle y,\lambda_n\rangle}e^{i\langle x,\lambda_n\rangle}
=\sum_{n=0}^{\infty}a_n e ^ {i\langle z, \lambda_n\rangle},
\ a_n\in {{\bf'}},
\end {equation}
(see~\cite {Ron1}).
Any series of the form~(\ref {f+}) is called {\it Diriáhlet series }.

By $\Gamma $ we always denote a convex closed cone
in $ {{\bf R}^m} $; by
$\widehat{\Gamma}$ we denote {\it the conjugate cone to $ \Gamma $:}
$$\widehat{\Gamma}=
\{t\in {{\bf R}^m}: \langle t, y\rangle \ge 0\ \forall y\in \Gamma \},$$
note that $ \widehat{\widehat{\Gamma}}=\Gamma.$ Also,
$\stackrel{\circ}{\Gamma}$ is the interior of a cone $ \Gamma. $

\begin {th-ud} \hspace {-0,5em}.
\label {th-ud8}
Let $f(x) $ be an almost periodic function by Bohr on $ {{\bf R}^m} $
with Fourier series
\begin{equation}
\label{f5}
\sum_{n=0}^{\infty}
a_n e^{i\langle x,\lambda_n\rangle},
\end{equation}
where all the exponents $\lambda_n$ belong
 to a cone $ \Gamma\subset {{\bf R}^m} $.
Then $f(x)$ continuously extends to the tube set
$T_{\widehat{\Gamma}}$ as an
 almost periodic by Bohr function
 $F (z) $ with Fourier series~(\ref{f+}).
The function $F(z) $ is holomorphic on the
interior $T_{\widehat{\Gamma}}$, and for any
$\Gamma'\subset\subset\widehat{\Gamma}$ uniformly w.r.t.
$z\in T_{\Gamma'}$
\begin{equation}
\label{lim}
\lim\limits_{\|y\|\rightarrow\infty}F(z)=a_0,
\end {equation}
where $a_0 $ is the Fourier coefficient corresponding to the
 exponent $\lambda_0=0$
(if \  $0\notin~sp~\,f $, put $a_0=0. $)
If $sp\,f\subset\,\stackrel{\circ}{\Gamma}$, then (\ref {lim}) is true
uniformly w.r.t. $z\in T_{\widehat \Gamma} $.
\end {th-ud}

Here the inclusion $ \Gamma'\subset\subset\widehat\Gamma $ means
that the intersection of
$\Gamma'$ with the unit
 sphere is contained in the interior of the intersection of $\widehat{\Gamma}$ with this
sphere.

To prove this theorem, we use the following lemmas.

\begin{lem-ud}\hspace{-0,5em}.
\label{l0}
Suppose that a plurisubharmonic function $ \varphi (z) $ on ${\bf C}^m $
is bounded from above on a set $T_{K} $, where $K\subset{\bf R}^m $ is
a convex set.
 Then the function
$$\psi(y)=\sup\limits_{x\in {{\bf R}^m}}\varphi(x+iy)$$ is convex on $K. $
\end {lem-ud}

{ \it Proof of Lemma~{\ref {l0}}.}
Fix $y_1, \ y_2\in K $.
The plurisubharmonic on ${\bf C}^m $ function
$$\varphi_1(z)=\varphi(z)-\frac{\psi(y_2)-\psi(y_1)}{\|y_2-y_1\|^ 2}\
\langle {\rm Im \, z}, y_2-y_1\rangle $$
is bounded from above on the set
$T_{[y_1, y_2]} $.
Therefore the subharmonic function $\varphi_2(w)=\varphi_1((y_2-y_1)w+iy_1)$
is bounded on the strip
$\{w=u+iv: u\in{\bf R}, \ 0\le v\le 1 \} $.
Hence, the value of $\varphi_1$ at
any point of this strip does not exceed
$$
\max\{\sup\limits_{u\in {\bf R}}\varphi_2(u),
\sup_{u\in {\bf R}}\varphi_2 (u+i)\}
\le
\max\{\sup\limits_{x\in {\bf R}^m}\varphi_1(x+iy_1),
\sup_{x\in {\bf R}^m}\varphi_1(x+iy_2)\}. $$
Therefore, for any $z=x+iy\in T_{[y_1, y_2]}$,
$$
\varphi_1(z)\le\max\{\psi(y_1)-\frac{\psi(y_2)-\psi(y_1)}{\|y_2-y_1\|^2}\
\langle y_1, y_2-y_1\rangle,
\psi(y_2)-\frac{\psi(y_2)-\psi(y_1)}{\|y_2-y_1\|^2}\
\langle y_2,y_2-y_1\rangle\}$$
$$=\frac{\|y_2\|^2-\langle y_1, y_2\rangle}{\|y_2-y_1\|^2}\ \psi(y_1)+
\frac{\|y_1\|^2-\langle y_1, y_2\rangle}{\|y_2-y_1\|^2}\ \psi(y_2).$$
Hence for any $y\in [y_1,y_2]$
we have
$$
\psi(y)\le
\frac{\|y_2\|^2-\langle y_1, y_2\rangle-\langle y, y_2-y_1\rangle}
{ \|y_2-y_1\|^2} \ \psi(y_1)+
\frac{\|y_1\|^2-\langle y_1, y_2\rangle +\langle y, y_2-y_1\rangle}
{\|y_2-y_1\|^2}\ \psi(y_2).$$
If $y =\lambda y_1 + (1-\lambda) y_2 $, $ \lambda\in (0,1) $, then
we obtain the inequality
$$
\psi(\lambda y_1+(1-\lambda) y_2)\le\lambda\psi(y_1)+(1-\lambda)\psi(y_2).$$
Therefore, the function $\psi(y)$ is convex on $K$.

\vspace {-1em}
\qed

\begin {lem-ud} \hspace {-0,5em}.
\label {l4}
Let $\psi(y)$ be a convex bounded
function on a cone $\Gamma$. Then
$\psi(y)\le\psi(0)$ for all $y\in \Gamma $.
\end {lem-ud}

{ \it Proof of Lemma~{\ref {l4}}.}
Since $\psi(y)$ is convex, we have
$$\psi(y)\le
\left(1-\frac 1t
\right)
\psi(0)+\frac 1t
\psi(ty),\ t>1.
$$
Taking $t\rightarrow \infty $,
we obtain $\psi(y)\le\psi(0).$ \\

\vspace {-1em}
\qed

{ \it Proof of Theorem~{\ref{th-ud8}}.}
Let $\sigma_q (x), \ q=0,1,2,\ldots $ be the Bochner-Feyer sums
 for the series~(\ref {f5}).
Obviously, these functions are also defined for $z\in {\bf C}^m $.
Assume that
$$\varphi_{q,l}(z)=\log(|\sigma_q(z)-\sigma_l(z)|).$$
For any fixed $q$ and $l$, $q>l$, the function $ \varphi_{q,l}(z)$
is plurisubharmonic on $ {\bf C}^m $. Moreover, for
$z\in T_{\widehat{\Gamma}}$ we have
$\langle y, \lambda_n\rangle \ge 0 $ and
$$
|\sigma_q(z) -\sigma_l(z)|\le |\sigma_q (z)|+|\sigma_l(z)|\le$$
$$\sum\limits_{n=0}^{q-1}|a_n|e^{-\langle y,\lambda_n\rangle}+
\sum\limits_{n=0}^{l-1}|a_n|e^{-\langle y,\lambda_n\rangle}\le
2\sum\limits_{n=0}^{q-1}|a_n|.$$
Consider the function
$$\psi_{q,l}(y)=
\sup_{x\in{{\bf R}^m}}\log(|\sigma_q(z)-\sigma_l(z)|).$$
Using lemma~\ref {l0}, we obtain that
$\psi_{q,l}(y)$ is convex in $ \widehat {\Gamma}. $
Therefore, by lemma~\ref {l4}, we have
\begin{equation}
\label{sup}
\sup\limits_{z\in T_{\widehat{\Gamma}}}(|\sigma_q(z)-\sigma_l(z)|)\le
\sup\limits_{x\in {{\bf R}^m}}(|\sigma_q(x)-\sigma_l(x)|).
\end {equation}
Further, the function $f(x)$ is almost periodic,
therefore the Bochner-Feyer sums converge
uniformly on $ {\bf R}^m$,
and for
$q, l\ge N(\varepsilon) $
$$\sup\limits_{x\in{{\bf R}^m}}(|\sigma_q(x)-\sigma_l(x)|)\le\varepsilon.$$
Hence,
$\sup\limits_{x\in {{\bf R}^m}}(|\sigma_q(z)-\sigma_l(z)|)\le\varepsilon$
for all $z\in T _ {\widehat {\Gamma}} $, $q, l\ge N (\varepsilon) $.

Thus the Bochner-Feyer sums uniformly converge  on $T_{\widehat{\Gamma}}$,
 their limit is an almost periodic
function by Bohr and holomorphic on the interior of $T_{\widehat{\Gamma}}$
with the Diriáhlet series~(\ref{f+}).

Further, passing to the limit in (\ref{sup}) as $q\rightarrow \infty $, we get
$$\sup_{z\in T_{\widehat{\Gamma}}}(|F(z)-\sigma_l(z)|)\le
\sup_{x\in {{\bf R}^m}}(|f(x)-\sigma_l(x)|).$$
Choose $l $ such that the right hand side of this inequality is less than
 $ \varepsilon $. We have for $\Gamma'\subset\subset{\widehat\Gamma}$
$$
\sup_{z\in T_{\Gamma'}}|F(z)-a_0|\le
\sup_{z\in T_{\Gamma'}}
|F(z)-\sigma_l(z)|+\sup_{z\in T_{\Gamma'}}|\sigma_l(z)-a_0|
\le\varepsilon+\sup\limits_{z\in T_{\Gamma'}}|\sigma_l(z)-a_0|.
$$
Note that for any fixed
$\lambda_n\in\Gamma\setminus\{0\}$
the value $ \langle y,\lambda_n\rangle$
tends to $+\infty $ as $\|y\|\rightarrow\infty$,
$y\in\Gamma'$, therefore, we have
$$|\sigma_l(z)-a_0|=|\sum\limits_{j=1}^{l-1}k_j^l
a_je^{i\langle x,\lambda_j\rangle}
e^{-\langle y,\lambda_j\rangle}|\le
\sum\limits_{j=1}^{l-1}|a_j|e^{-\langle y,\lambda_j\rangle}
\rightarrow 0$$
as $\|y\|\rightarrow\infty$ on $\Gamma'$.
Hence, uniformly w.r.t. $z\in T_{\Gamma'}$
\begin{equation}
\label{mod}
\overline{\lim\limits_{\|y\|\rightarrow\infty}}
|F(z)-a_0|\le\varepsilon.
\end {equation}
This is true for arbitrary $\varepsilon > 0 $, then (\ref{lim})
follows.

If $sp\,f\subset\ \stackrel{\circ}{\Gamma}$, then for any
 $ \lambda_n\in sp\,f$,\ \
$ \langle y, \lambda_n\rangle\rightarrow+\infty $ as
$\|y\|\rightarrow\infty$ uniformly w.r.t. $y\in \widehat{\Gamma} $,
therefore (\ref{mod}) is true uniformly w.r.t.
 $z\in T_{\widehat{\Gamma}}$, and (\ref {lim}) is also true.
The theorem has been proved.

\qed

\begin {th-ud} \hspace {-0,5em}.
\label {th-ud9}
Let $f(x)$ be an almost periodic function by Stepanoff on $ {{\bf R}^m}$
with the Fourier series~(\ref {f5}). Let
all the exponents $ \lambda_n $ belong to a cone $\Gamma\subset{{\bf R}^m}$.
Then there exists an almost periodic by Stepanoff function $F(z)$ in the
tube set $T_{\widehat{\Gamma}}$ with the Fourier series~(\ref{f+})
 such that $F(x)=f(x)$. The function $F(z)$
is holomorphic almost periodic by Bohr
on any domain $T_{\widehat{\Gamma}+b}$,
$b\in \stackrel{\circ}{\widehat{\Gamma}}$.
Besides, for any cone $ \Gamma'\subset\subset{\widehat\Gamma}$
we have uniformly w.r.t. $z\in T_{\Gamma'}$
\begin{equation}
\label{a_0}
\lim\limits_{{\|y\|\rightarrow\infty}}F(z)=a_0,
\end {equation}
where $a_0 $ is the Fourier coefficient for the exponent $\lambda=0. $
If $sp\,f\subset\stackrel{\circ}{\Gamma}$, then (\ref{a_0})
is true uniformly w.r.t. $z\in T_{\widehat{\Gamma}+b} $ for any
$b\in\stackrel{\circ}{\widehat\Gamma}$.
\end {th-ud}

{ \it Proof.}
To prove the first part of the theorem, we need to replace
$\varphi_{q,l}(z)$ by
$$\widehat{\varphi}_{q,l}(z)=
\log\left(
\int\limits_{[0,1]^m}^{}|\sigma_q(z+u)-\sigma_l(z+u)|^p du
\right)^\frac 1p.$$
Arguing as in the proof of theorem~\ref{th-ud8}, we obtain that the
 Bochner-Feyer sums $ \sigma_q (z) $
converge in the Stepanoff metric uniformly w.r.t.
$z\in T_{\widehat {\Gamma}}$ to
an almost periodic function by Stepanoff $F (z) $ with Fourier
series~(\ref{f+}).

Let $b\in \stackrel{\circ}{\widehat\Gamma}$.
The module of the function $\sigma_q(z)-\sigma_l(z)$
is estimated from above by
the mean value on the corresponding ball contained
 in $T_{\widehat{\Gamma}}$.
Using
the H\"older inequality, we have
$$\sup\limits_{x\in {{\bf R}^m}}|\sigma_q(x+bi)-\sigma_l(x+bi)|
\le C\sup\limits_{z\in T_{\widehat{\Gamma}}}
\left(
\int\limits_{[0,1]^m}{}|\sigma_q (z+u) -\sigma_l (z+u)|^p du
\right)^{\frac 1p}, $$
where the constant $C $  depends only on $b $ and $ \widehat{\Gamma}.$

Applying lemmas~\ref{l0} and \ref{l4} to the functions
$$
\tilde\psi_{q,l,b}(y)=\sup\limits_{x\in {\bf R}^m}\log|\sigma_q (z+bi)-\sigma_l(z+bi)|, $$
we get that the Bochner-Fourier sums converge uniformly on
$T_{\widehat {\Gamma}+b} $
to $F (z) $, thus $F (z) $ is holomorphic almost periodic by Bohr
in $T_{\widehat{\Gamma}+b}$ for any
$b\in\stackrel{\circ}{\widehat\Gamma}$.

Then the other statements of the theorem
 follow from theorem~\ref{th-ud8}.

\qed
\vspace {-1em}

Now we prove the inverse statements to theorems~\ref {th-ud8}
and~\ref {th-ud9}.

\begin {th-ud} \hspace {-0,5em}.
\label {th-ud11}
Suppose that an almost periodic by Bohr function $f(x)$
continuously extends to the interior of
$T_{\Gamma}$ as a holomorphic function $F(z)$. If $F (z) $ is bounded on any set
$T_{\Gamma'}$, $\Gamma'$ being a the cone in $ {\bf R}^m $, $\Gamma'\subset\subset\Gamma,$
then $F (z) $ is an almost periodic function by Bohr on $T_{\Gamma}$
and the spectrum of $F(z)$ is contained in $\widehat{\Gamma} $.
\end {th-ud}

{ \it Proof.}
Take
$ \lambda\notin \widehat {\Gamma} $. Then there exists
$y_0\in\ \stackrel{\circ}{\Gamma}$ such that $\langle y_0,\lambda\rangle <0$.

Choose a neighbourhood $U\subset\stackrel{\circ}{\Gamma}$ of $y_0$ such that
$\langle y,\lambda\rangle\le\frac 12 \langle y_0,\lambda\rangle$
for all $y\in U.$
Let $A $ be any nondegenerate operator in $ {\bf R}^m $ such that $A$
  maps all the vectors
$e_1=(1,0,\ldots,0),\ldots, e_m=(0,0,\ldots,1)$
into $U$.

The function $F(A\zeta)$
 is holomorphic and bounded on the set
$$\{\zeta=\xi+i\eta \in {\bf C}^m: \xi\in{\bf R}^m,\
\eta^j>0,j=1,\ldots,m \}$$
because
$A\{\eta : \eta^j\ge 0\}\subset\Gamma$.

If for each coordinates $\zeta^1,\ldots,\zeta^m$ we change the
integration over the segments $-N\le\xi^j\le N$, $\eta^j=0$
to the integration over the half-circles $\zeta^j=N e^{i\theta^j}$,
$0\le\theta^j\le\pi$, $j=1,\ldots,m$
we obtain the equality
\begin{equation}
\label{prod}
\left(\frac{1}{2N}\right)^m
\int\limits_{[-N,N]^m}^{}
F(A\xi) e^{-i \langle A\xi,\lambda\rangle}d\xi=
\left(\frac{i}{2}\right)^m\int\limits_{[0,\pi]^m}^{}
F(ANe^{i\theta})\prod\limits_{j=1}^{m}e^{i\theta^j-iNe^{i\theta^j}
\langle Ae_j,\lambda\rangle}d\theta,
\end{equation}
where $\theta=(\theta^1,\ldots,\theta^m),$
$e^{i\theta}=(e^{i\theta^1},\ldots,e^{i\theta^m})$.

Since $\langle Ae_j,\lambda\rangle <0$ for $j=1,\ldots,m$, we see
that the integrand in the right-hand side of (\ref{prod}) is uniformly
bounded for all $N>1$. By Lebesgue theorem (\ref{prod}) tends to zero as
$N\rightarrow\infty$.

Thus
\begin{equation}
\label{=0}
{\lim\limits_{N\rightarrow\infty}}
\frac {1}{(2N)^m} \int\limits_{A ([-N, N]^m)}
f (x) e^{-i\langle x, \lambda\rangle} dx=0.
\end {equation}
Cover the set $A ([-N^2, N^2] ^m) $ by cubes $L_j=x'_j+[-N, N]^m $
such that the interiors of these cubes are not intersected.
We may assume that the number of the cubes intersecting the
boundary of the set $A([-N^2,N^2]^m)$ is
$O(N^{m-1}) $ as $N\rightarrow\infty $.
Taking into account boundedness of the function
$f(x)e^{-i\langle x, \lambda\rangle}$ on ${{\bf R}^m}$ and
 equality~(\ref{=0}),
we have
\begin{equation}
\label{N^{2m}}
\frac {1}{(2N)^{2m}}
\int\limits_{\cup L_j}f (x)e^{-i\langle x, \lambda\rangle} dx =
\frac{1}{(2N)^{2m}}
\left(\int\limits_{A ([-N^2, N^2]^m)}
f (x) e^{-i\langle x, \lambda\rangle} dx+O(N^{2m-1})\right)=o(1)
\end{equation}
Since~(\ref{K})
we see that uniformly w.r.t. $x'\in {\bf R}^m $
as $N\rightarrow\infty$
\begin{equation}
\label{intx'}
\frac {1}{(2N)^{m}}
\int\limits_{x'+[-N, N]^m}f(x)e^{-i\langle x, \lambda\rangle} dx
=a(\lambda, f)+o(1).
\end{equation}

On the other hand, the number of the cubes $L_j $ equals $O(N^{m}) $ as
$N\rightarrow\infty $, then the equality $a(\lambda, f)=0$  follows
 from (\ref{N^{2m}}) and (\ref{intx'}).
This yields the inclusion
$sp\,f\subset\widehat{\Gamma}. $
Using theorem~\ref {th-ud8} we complete the proof of our theorem.

\qed
\vspace {-1em}

\begin {th-ud} \hspace {-0,5em}.
\label {th-ud12}
If $F (z) $ is bounded on each set $T_{\Gamma'} $,
$\Gamma'\subset\subset\Gamma$, and the
nontangential limit value of $F(z)$ as $y\rightarrow 0$
is an almost periodic function by Stepanoff on
${\bf R}^m$, then $F(z)$ extends to $T_{\Gamma} $
as an almost periodic function by
Stepanoff and the spectrum of $F(z)$ is contained in $\widehat{\Gamma} $.
\end{th-ud}

{ \it Proof.}
The proof of this theorem is the same as of theorem~\ref{th-ud11},
 but we have to use theorem~\ref{th-ud9} instead of theorem~\ref {th-ud8}.

\vspace {-1em}
\qed\\

To formulate further results we need the
concept of $P$-indicator.
( See, for example,~\cite {Ronkin}, p. 275.)\\
{\bf Definition.} {\it
$P$-indicator of an entire function
$F(z)$ on ${\bf C}^m$ is the function
$$h_F(y) =\sup\limits_{x\in {{\bf R}^m}}
\overline{\lim\limits_{r\rightarrow\infty}}
\frac 1r \log|F(x+iry)|.$$}

\begin {th-ud} \hspace {-0,5em}.
\label {th-ud10}
(For $m=1 $ see~\cite{Levin}, \cite{Levitan}.)
Let $f (x) $, $x\in {\bf R} ^m $ be an almost periodic function by Stepanoff
with the Fourier series~(\ref {f5}), and let $\|\lambda_n\|\le C<\infty$ for all~$n.$
Then $f(x)$ extends to ${\bf C}^m$ as an entire function
$F(z)$ of exponential type, which is almost periodic by Bohr
on any tube domain in $ {{\bf C}^m} $ with bounded base; $F (z) $
has the Dirichlet series~(\ref{f+}), and
 $P$-indicator $h_F(y)$ satisfies the equation
$h_F(y)=~H_{sp\, f}(-y)$, where
$H_{sp \, f}(\mu):=\sup\limits_{x\in sp\, f} \langle x, \mu\rangle $
is the support function of the set $sp\,f$.
\end {th-ud}

{ \it Proof.}
Take $ \mu\in{{\bf R}^m}$ such that $\|\mu\|=1$. Put
$$f_{\mu}(x)=f(x)e^{-i[H_{sp\,f}(\mu)+\varepsilon]\langle x,\mu\rangle}.$$
The Fourier series
$\sum\limits_{n=0}^{\infty}a_ne^{i\langle x,\lambda_n-(H_{sp\,f}(\mu)+
\varepsilon)\mu\rangle}$ corresponds to the function $f_{\mu}(x)$,
 hence $$sp\,f_{\mu}\subset
\{x\in{{\bf R}^m}:\langle x,\mu\rangle\ \le -\varepsilon\}.$$
Since $sp \, f_{\mu} $ is bounded, we obtain for some~$\delta>0 $
$$
sp\,f_{\mu}(x)\subset\Gamma_{\delta,-\mu}=
\{\lambda\in {\bf R}^m: \langle\lambda, -\mu\rangle\ge\delta\|\lambda\|\}.$$
Theorem~\ref{th-ud9} yields that $f_{\mu}(x)$ extends
to the interior of the domain
$T_{\widehat{\Gamma}_{\delta, -\mu}},$ where
$$\widehat{\Gamma}_{\delta,-\mu}=
\{y:\langle y,-\mu\rangle\ge\sqrt{1-{\delta^2}}\|y\|\}$$ is
the conjugate cone to ${\Gamma}_{\delta,-\mu},$
as an almost periodic function by Bohr $F_{\mu}(z)$.
This function is holomorphic on any domain
$T_{\widehat{\Gamma}_{\delta,-\mu}+b}$,
$b\in\stackrel{\circ}{\widehat{\Gamma}}_{\delta,-\mu}$
 with the Dirichlet series
$$
\sum\limits_{n=0}^{\infty}a_n
e^{i\langle z,\lambda_n-[H_{sp\,f}(\mu)+\varepsilon]\mu\rangle},$$
and $F_{\mu}(z)\rightarrow 0 $
as $\|y\|\rightarrow\infty $
uniformly w.r.t.
$z\in T_{\Gamma'}$ for any cone
$\Gamma'~\subset\subset~\widehat{\Gamma}_{\delta,-\mu}$.
Using (\ref{K}) we get
\begin{equation}
\label{for15}
\left|
a_ne^{-\langle y, \lambda_n-[H_{sp\,f}(\mu)+\varepsilon]\mu\rangle}
\right|\le
\sup\limits_{x\in{\bf R}^m}|F_\mu(x+iy)|,\ y\in\Gamma'.
\end{equation}
Put
$$ F(z):=F_{\mu}(z)e^{i[H_{sp\,f}(\mu)+\varepsilon]\langle z,\mu\rangle}.$$
$F(z)$ is almost periodic on $T_{\Gamma'} $ with
 Dirichlet series~(\ref{f+}).
Therefore it follows from (\ref{for15}) that
\begin{equation}
\label{an}
|a_n|\le\sup\limits_{x\in{{\bf R}^m}}
|F(x+iy)|e^{\langle y, \lambda_n\rangle}.
\end{equation}
On the other hand, the function $F_{\mu}(z)$ is bounded on
$T_{\Gamma'}$, hence
\begin{equation}
\label{f6}
|F(z)|
\le C(\Gamma')e^{-[H_{sp\,f}(\mu)+\varepsilon]\langle y,\mu\rangle},
\ z\in T_{\Gamma'}
\end{equation}

Cover the space ${{\bf R}^m} $ by the interiors of a finite number of cones
 $\Gamma_1',\ldots,\Gamma_N'$. There exist
 holomorphic on the interior of $\Gamma_k'$
almost periodic functions $F_k (z) $, $k=1, \ldots, N $ with
 identical Dirichlet series~(\ref{f+}).
 Using the uniqueness theorem, we obtain that these functions coincide
 on the intersections of the cones and thus define a holomorphic
function $F (z) $ on ${{\bf C}^m}\setminus{{\bf R}^m} $.
The Bochner-Feyer sums for $F(z)$
converge to this function  uniformly on any set
$$\{z=x+iy:x\in {\bf R}^m,\ \|y\|=r>0\}.$$
Hence, these sums converge on the tube domain $T_{\{\|y\|<r\}}$.
Thus
$F(z)$ extends to $ {\bf C}^m $ as the holomorphic function, which is
 almost periodic on any tube set with a bounded base.
Owing to the uniqueness of expansion into Fourier series, we have $F(x)=f(x)$.

Let us prove that $h_F(y)=H_{sp \, f}(-y)$.
From inequality~(\ref {f6}) with $ \mu =-y $ it follows that
$$h_F(y)\le\overline{\lim\limits_{r\rightarrow\infty}} \frac 1r
[ H_{sp \, f} (-y) + \varepsilon] \langle ry, y\rangle=
H_{sp \, f} (-y) + \varepsilon. $$
The functions $h_F(y)$ and
$H_{sp \, f} (y) $
are positively homogenious, hence
the inequality
$$ h_F (y) \le H_{sp \, f} (-y) $$ is true for all $y\in {\bf R}^m. $

Further, fix $x, y\in {\bf R}^m $.
The holomorphic on ${{\bf C}} $ function
$ \varphi (w) =F (x+wy) $ is  bounded
on the axis
${\rm Im}\,w=0$.
Then the estimate
$$ | \varphi (w) | \le Ce ^ {a | {\rm Im} \, w |} $$
for some $a > 0 $ and all $w\in {\bf C} $
follows from~(\ref {f6}).
Using the definition of $P$-indicator, we get
$$\overline {\lim\limits_{v\rightarrow + \infty}}
\frac1v\ \log|\varphi(iv)|\le h_F(y). $$
Therefore the function
$\varphi(w)e^{i(h_F(y) + \varepsilon) w} $ is bounded on
the positive part of the imaginary axis.
Applying the Fragmen-Lindelof principle to the quadrants
${\rm Re}\,w\ge 0,\ {\rm Im}\,w\ge 0 $ and
${\rm Re} \, w\le 0,\ {\rm Im} \, w\ge 0 $,
we get boundedness of this function on the upper half-plane.
Applying the Fragmen-Lindelof principle to the half-plane
 ${\rm Im}\,w\ge 0 $, we get the inequality
$$
|\varphi (w) | \le
\left(
\sup\limits_{{\rm Im} \, w=0} | \varphi (w) |
\right)
e^{h_F(y){\rm Im}\,w}
\eqno({\rm Im}\,w>0).$$
Hence, for all $z\in {{\bf C}^m} $, we have
$$ |F (z) |
\le \sup\limits_{x\in {{\bf R}^m}}|F (x) |e^{h_F (y)}.$$
Now using formula~(\ref{an}) for coefficients of the Dirichlet series
of the function $F (z)$,
we get the estimate
\begin{equation}
\label{f7}
|a_n|\le \sup\limits_{x\in {{\bf R}^m}}
|f (x) |e^{h_F (y) + \langle y, \lambda_n\rangle}.
\end{equation}
Suppose
$ \langle y_0, \lambda_n\rangle +h_F (y_0)~<~0 $ for some
$y_0\in {\bf R}^m $. Put $y=ty_0 $ in (\ref {f7}) and
let $t\rightarrow \infty $. We obtain
$a_n=0 $. This is impossible because $ \lambda_n\in sp \, f $.

Thus for all $y\in {{\bf R}^m} $ and
$ \lambda_n\in sp\,f $ we have
$h_F (y) + \langle y, \lambda_n\rangle \ge 0 $, hence
$$ H_{spf}(-y) = \sup\limits _ {\lambda_n\in sp \, f}
\langle -y, \lambda_n\rangle \le h_F (y). $$
This completes the proof of the theorem.

\vspace {-1em}
\qed

The following theorem  is inverse
 to the previous one.
\begin{th-ud}\hspace{-0,5em}.
\label{th-ud13}
(For $m=1 $ see~\cite {Levin}, \cite {Levitan}.)
Let $F (z) $ be an entire function on $ {{\bf C}^m} $,
$|F(z)|\le Ce^{b\|z\|},$
let $F (x) $, $x\in{\bf R}^m$ be an almost periodic function by
Stepanoff with the Fourier series~(\ref {f5}).
Then $F (z) $ is an almost periodic function by Bohr on any tube domain
$T_{D}\subset{{\bf C}^m}$ with the bounded base, $F (z) $
has the Dirichlet series~(\ref{f+}),
 and $sp\,F~\subset~\{\lambda:\|\lambda\|\le b\} $.
\end{th-ud}

{\it Proof.}
It follows from theorem~\ref {th-ud10}, that
 it suffices to prove the inclusion
$$ sp\,F\subset\{\lambda:\|\lambda\|\le b \}.$$
Let the function $F(x)$ be bounded on $ {\bf R}^m $.
Arguing as in theorem~\ref{th-ud10},
 we see that for all $z\in {\bf C}^m $
$$ |F (z) | \le\sup\limits_{x\in {\bf R}^m} |F (x) |e^{h_F (y)}, $$
where $h_F (y) $ is $P$-indicator for $F(z).$
Further, for all $x\in {\bf R}^m $ we have
$$
h_F (y) = \sup_{x\in {\bf R}^m}
\overline{\lim\limits_{r\rightarrow\infty}}
\frac 1r \log|F(x+iry)|\le\sup_{x\in {\bf R}^m} \overline{\lim
\limits_{r\rightarrow\infty}}
\frac 1r (\log C+b\| x+iry\|)\le b\|y\|, $$
therefore for all $z\in {\bf C}^m $
$$ |F(z)|\le Ce^{b \|y\|}. $$
Take $ \varepsilon > 0 $, $ \mu\in {{\bf R}^m}, \, \|\mu\| = 1. $
Consider the function
$$ F_{\mu}(z)=F(z)e^{-i\langle z, \mu\rangle (b+\varepsilon)}. $$
Since $ |F_{\mu}(z)|\le Ce^{b\|y\|+(b+\varepsilon)\langle y,\mu\rangle}$
 uniformly w.r.t. $x\in {\bf R}^m $, then
$F_{\mu}(z)$ is uniformly bounded for
$z\in T_{\Gamma_{-\mu}}$, where $\Gamma_{-\mu}$ is the cone
$\{y:\langle y,-\mu\rangle\ge(1-\frac{\varepsilon}{b+\varepsilon})\|y\|\}.$
Using theorem \ref {th-ud12}, we obtain that the spectrum
$F_\mu$ is contained in
$\widehat{\Gamma}_{-\mu}$ and
$$sp\,F=sp\,F_{\mu}+(b+\varepsilon)\mu\subset\widehat{\Gamma}_{-\mu}+
(b+\varepsilon)\mu.$$

Finally, using the inclusion
$$\bigcap_{\mu:\,\|\mu\|=1}
(\widehat {\Gamma}_{-\mu}+(b+\varepsilon)\mu)
\subset\{\lambda:\|\lambda\|\le b+\varepsilon\}$$
and the arbitrarity of choice of $\varepsilon$ we get the assertion
of the theorem in the case of bounded on $ {\bf R}^m $ function $F(z)$.

Now let the function $F(z)$ be unbounded on ${\bf R}^m. $
Put for some $N> 0 $
$$
g(z)=\frac{1}{N^m}\int\limits_{[0,N]^m}^{}F(z+t)dt.$$
The function $g(z)$ satisfies the estimate on ${\bf C}^m $
$$
|g(z)|\le Ce^{bmN}e^{b\|z\|}.$$
As in the case $m=1$ (see ~\cite {Levitan}), we can prove
that
$g(x)$ is an almost periodic function by Bohr and is bounded
on ${\bf R}^m $.
The function $g (x) $ has
the Fourier series
$$
\sum\limits_{n=0}^{\infty}
a_n\frac{e^{i\lambda_n^1N}-1}{N\lambda_n^1}
\ldots
\frac{e^{i\lambda_n^mN}-1}{N\lambda_n^m}e^{i\langle x, \lambda_n\rangle},$$
where $\lambda_n^j$ are coordinates of the vector $\lambda_n$
(if $\lambda_n^j=0$,
the corresponding multiplier should be replaced by 1).

Using countability of $sp\,F $, we can choose $N$ in such a way that
none of the numbers
$\lambda_n^jN$ coincides with $2\pi k$, $k\in {\bf Z}\setminus\{0\}.$
In this case $sp\,g=sp\,F$.
Applying the proved above statement to the function $g (z) $, we
obtain the inclusion
$$ sp\,F\subset\{\lambda:\|\lambda\|\le b\}.$$
\qed

\renewcommand {\refname} {\normalsize \rm 
\end {document}